\documentclass[twocolumn]{autart}
\usepackage{amsmath}
\usepackage{amssymb,graphicx}
\usepackage{amsfonts}

\usepackage{natbib}

\newtheorem{Prop}{Proposition}

\newtheorem{Cor}{Corollary}
\newtheorem{rk}{Remark}

\newcommand{\wmin}{\omega_{*}}
\newcommand{\wmax}{\omega^{*}}

\newcommand{\SSS}{{\mathbb  S}}

\begin{document}
\bibliographystyle{elsarticle-harv}

\begin{frontmatter}
\runtitle{~}  

\title{Stabilization of an arbitrary  profile \\ for an
ensemble of half-spin systems\thanksref{footnoteinfo}} 

\thanks[footnoteinfo]{Corresponding author P. Rouchon Tel. +33 1 40 51 91 15.
Fax +33 1 40 51 91 65.}

\author[ENS]{Karine Beauchard}\ead{Karine.Beauchard@cmla.ens-cachan.fr},
\author[USP]{Paulo S\'ergio Pereira da Silva}\ead{paulo@lac.usp.br},    
\author[MINES]{Pierre Rouchon}\ead{pierre.rouchon@mines-paristech.fr}              

\address[ENS]{CMLA, ENS Cachan, CNRS, UniverSud, 61, avenue du Pr\'{e}sident Wilson,
F-94230 Cachan, FRANCE.}

\address[USP]{University of S\~ao Paulo, Escola Polit\' ecnica --
PTC --- Av. Luciano Gualberto trav. 03, 158, 05508-900 -- S\~ao
Paulo -- SP BRAZIL}

\address[MINES] {Mines ParisTech, Centre Automatique et Syst\`{e}mes,
Unit\'{e} Math\'{e}matiques et Syst\`{e}mes,
 60 Bd Saint-Michel, 75272 Paris cedex 06, FRANCE}

\begin{keyword}                           
nonlinear systems, Lyapunov stabilization, LaSalle invariance, quantum systems, Bloch equations,   ensemble controllability, infinite dimensional system.
\end{keyword}                             

\begin{abstract}                          
We consider the feedback stabilization of a variable profile for an ensemble of non interacting half spins described by the Bloch equations. We propose an explicit feedback law that stabilizes asymptotically the system around a given arbitrary  target profile. The convergence  proof is done when the target profile  is entirely in the south hemisphere or in the north hemisphere of the Bloch sphere. The convergence holds for initial conditions in a $H^1$ neighborhood of this target profile. This convergence is shown for the weak $H^1$ topology. The proof relies on an adaptation of the LaSalle invariance principle to infinite dimensional systems. Numerical simulations illustrate the efficiency of these feedback laws, even for initial conditions far from the target profile.
\end{abstract}
\end{frontmatter}

\section{Introduction}
Ensemble controllability as introduced  in~\cite{li-khaneja:ieee09} is an  interesting control theoretic notion  well adapted to  nuclear magnetic resonance (NMR) systems (see, e.g., \cite{li-khaneja:pra06} and the reference herein). In~\cite{beauchard-et-al:JMP2010} some controllability issues of such NMR systems are investigated using open-loop controls involving  Dirac-combs. In \cite{beauchard-et-al:auto11}  such open-loop  Dirac-combs  are combined with  Lyapunov stabilizing feedback to ensure closed-loop convergence towards a target profile  that is  one of the two  steady-states, the south and  north poles of the Bloch sphere. In this note, we extend this Lyapunov design  to arbitrary target profiles and prove its local convergence for  weak $H^1$ topology when the target profile  lies entirely in the south hemisphere or in the north hemisphere.

We consider an ensemble of non interacting half-spins in a static
field $(0,0,B_0)^t$ in $\mathbb{R}^3$, subject to a transverse radio
frequency field $(\tilde{u}_1(t),\tilde{u}_2(t),0)^t$ in $\mathbb{R}^3$
(the control input). The ensemble of half-spins is described  by the
magnetization vector $M\in\mathbb{R}^3$ depending on time $t$ but
also on the Larmor frequency $\omega=-\gamma B_0$ ($\gamma$ is the
gyromagnetic ratio). It obeys to  the Bloch equation:
\begin{equation}
 \label{dyn:eq1}
  \frac{\partial M}{\partial t}(t, \omega) =
(\tilde u_1(t) e_1 + \tilde u_2(t) e_2 + \omega e_3) \wedge
  M(t, \omega),
\end{equation}
where $-\infty < \wmin < \wmax < +\infty$, $\omega\in(\wmin,\wmax)$,
$(e_1,e_2,e_3)$ is the canonical basis of $\mathbb{R}^3$, $\wedge$
denotes the wedge  product on $\mathbb{R}^3$. The equation
(\ref{dyn:eq1}) is an infinite dimensional bilinear control system.
The state is the $\omega$-profile  $M$, where, for every $\omega \in (\wmin,\wmax)$,
$M(t,\omega)\in  \mathbb{S}^{2}$ (the unit sphere of $\mathbb{R}^3$).
The two control inputs  $\tilde{u}_1$ and $\tilde{u}_2$ are real valued.

We propose here  a first answer to the  local stabilization  of an arbitrary  profile:
given an arbitrary target profile
$M_f:(\omega_*,\omega^*) \rightarrow \SSS^2$, define
an explicit control law $(\tilde{u}_1(t,M),\tilde{u}_2(t,M))$,
a neighborhood $U$ of $M_f$ (in some space of functions to be determined),
a diverging sequence of times $(t_n)_{n \in \mathbb{N}}$,
such that, for every initial condition $M^0 \in U$, the solution of
the closed loop system is uniquely defined and satisfies
$$\lim\limits_{n \rightarrow + \infty}
\| M(t_n,.) - M_f(.) \|_{L^\infty(\omega_*,\omega^*)} =0.$$
In this note, the Lyapunov feedback  proposed in~\cite{beauchard-et-al:auto11} is adapted to provide a  constructive answer to this question. Section~\ref{H1:sec}  is devoted to control design and closed-loop
simulations. In section~\ref{main:sec} we state and
prove the main convergence result, theorem~\ref{MainThm}.

\section{Lyapunov $H^1$ approach}{\label{H1:sec}
\subsection{Some preliminaries}

Let us recall the concept of a solution for (\ref{dyn:eq1}) when the control input  $u$ contains Dirac distributions.
When $\tilde{u}_1, \tilde{u}_2 \in L^1_{loc}(\mathbb{R})$,
then, for every initial condition $M_0 \in L^2((\omega_*,\omega^*),\mathbb{R}^3)$,
the equation (\ref{dyn:eq1}) has a unique weak solution
$M \in C^0([0,+\infty),L^2((\wmin,\wmax),\mathbb{R}^3))$. Denote by $\delta(t-a)$ the Dirac distribution located at $t=a$.
When $\tilde{u}_1=\alpha \delta(t-a)+u_1^{\sharp}$ and
$\tilde{u}_2=u_2^{\sharp}$ where $u_j^{\sharp} \in
L^1_{loc}(\mathbb{R})$, $\alpha>0$ and $a \in (0,+\infty)$, then the
solution is the classical solution on $[0,a)$ and $(a,+\infty)$, it
is discontinuous at the time $t=a$, with an explicit discontinuity
given by an instantaneous rotation of angle $\alpha$ around the axis
$\mathbb{R}e_1$
$$M(a^+,\omega)=\left( \begin{array}{ccc}
1 & 0            &  0           \\
0 & \cos(\alpha) & -\sin(\alpha)\\
0 & \sin(\alpha) &  \cos(\alpha)
\end{array}\right) M(a^-,\omega).$$

The symbol $\|.\|$ (resp. $\langle . , . \rangle$) denotes the Euclidian norm (resp. scalar product)
on $\mathbb{R}^3$ and the associated operator norm on $\mathcal{M}_3(\mathbb{R})$.

\subsection{Transformation into a driftless system}

As in \cite{beauchard-et-al:auto11} we consider a control with an ``impulse-train'' structure
\begin{equation}\label{uv:eq}
\tilde {u}_1 = u_1 + \sum_{k=1}^{+\infty} \pi~\delta(t-kT), \quad
\tilde{u}_2 = (-1)^{\epsilon(t)} u_2
\end{equation}
where $\epsilon(t):=E(t/T)$, for some period $T>0$ and $E(\gamma)$ denotes the integer part of
the real number $\gamma$. The new controls $u_1, u_2$ belong to
$L^1_{loc}(\mathbb{R})$. Considering the change of variable
\begin{equation} \label{def:P(t)}
M_1(t,\omega) := P(t) M(t,\omega) \text{ where } P(t):=\left(
\begin{array}{ccc}
1 & 0           & 0          \\
0 & \epsilon(t) & 0          \\
0 & 0           & \epsilon(t)
\end{array} \right)
\end{equation}
one gets the following dynamics
\begin{equation} \label{Syst_ss_disp}
\frac{\partial M_1}{\partial t}(t,\omega) = [u_1(t) e_1 + u_2(t) e_2 +
\epsilon(t) \omega] \wedge M_1(t,\omega).
\end{equation}
The application of impulses at $t=kT$, by changing the sense of rotation of the null input solution,
is expected to reduce the dispersion in the closed loop system.
Since $M(t,\omega)=M_1(t,\omega)$ for every $t \in
[2kT,(2k+1)T]$, any convergence result on $M_1(t)$ when $t
\rightarrow + \infty$ provides a convergence result on $M$.

The first step of the control design consists in putting the system
(\ref{Syst_ss_disp}) in driftless form. The new function
$$M_2(t,\omega):=\exp[ \sigma(t) \omega S ] M_1(t,\omega)$$
where
\begin{equation} \label{def:S}
\sigma(t):=\int_0^t \epsilon(s) ds, \quad
S:=\left( \begin{array}{ccc}
0  & 1 & 0 \\
-1 & 0 & 0 \\
0  & 0 & 0
\end{array} \right),
\end{equation}
solves
\begin{equation} \label{equ:M_tilde}
\frac{\partial M_2}{\partial t}(t,\omega)= \sum_{i=1}^{2} u_i(t)
\Big[ \exp( \sigma(t) \omega S ) e_i \Big] \wedge M_2(t,\omega).
\end{equation}
Since $\sigma(2kT)=0, \forall k \in \mathbb{N}$, any convergence on
$M_2(t)$ when $t \rightarrow + \infty$ provides a convergence on
$M_1(2kT)$ when $k \rightarrow + \infty$.

\subsection{Transformation of the target profile}

The second step of the control design consists in transforming a
convergence to a variable profile $M_f$ into a convergence to the
constant profile $-e_3$, for which we developed tools in the
previous work \cite{beauchard-et-al:auto11}.
It relies on the following proposition.

\begin{Prop} \label{Def:R}
There exists $C>0$ such that, for all  $M_f\in H^1((\omega_*,\omega^*),\SSS^2)$,
there exists $R\in  H^1((\omega_*,\omega^*),SO_3(\mathbb{R}))$ satisfying
\begin{equation} \label{def:R}
R(\omega) M_f(\omega)=-e_3, \quad \forall \omega \in [\omega_*,\omega^*],
\end{equation}
\begin{equation} \label{R_borne}
\|R\|_{H^1} \leqslant C \|M_f\|_{H^1}.
\end{equation}
\end{Prop}

\noindent \textbf{Proof:} Let $M_f \in H^1((\omega_*,\omega^*),\SSS^2)$ and set
$f(\omega):=M_f'(\omega) \wedge M_f(\omega)$. Denote by $A(\omega)$ the skew-symmetric operator defined by
$\mathbb R^3 \ni M \mapsto f(\omega) \wedge M \in\mathbb R^3$. Consider the Cauchy problem
$$
\frac{d }{d \omega} R = R A(\omega)\text{ on } [\omega_*,\omega^*] \text{ with } R(\omega_*)=R_*
$$
where $R_*$ is any rotation sending $M_f(\omega_*)$ to $-e_3$: $R_* M_f(\omega_*)= - e_3$.  Since $\omega \mapsto A(\omega)$ is $L^2$ the solution $R$ is well defined,   unique  and belongs to $H^1((\omega_*,\omega^*),\SSS^2)$. Direct computations show that
$ \frac{d }{d \omega} (RM_f)  = 0$. Thus $R(\omega)M_f(\omega) \equiv  -e_3$. Moreover,
$\|R(\omega)\|=1$ and
$\|R'(\omega)\|=\|A(\omega)\|=\|M_f'(\omega)\|$ for all $\omega \in [\omega_*,\omega^*]$,
which proves (\ref{R_borne}). \hfill $\Box$.

Let us consider a target profile $M_f \in
H^1((\omega_*,\omega^*),\SSS^2)$. Take  $R\in  H^1((\omega_*,\omega^*),SO_3(\mathbb{R}))$ given by the above proposition. To
any solution $M_2$ of (\ref{equ:M_tilde}), we associate the function
\begin{equation} \label{def:N}
N(t,\omega):=R(\omega) M_2(t,\omega), \forall \omega \in
[-\omega_*,\omega^*].
\end{equation}
This function solves the equation
\begin{equation} \label{equ:N}
\frac{\partial N}{\partial t}(t,\omega)= \sum_{i=1}^{2} u_i(t) \Big[
F(t,\omega) e_i \Big] \wedge N(t,\omega)
\end{equation}
where
\begin{equation} \label{def:F}
F(t,\omega):=R(\omega) \exp( \sigma(t) \omega S ).
\end{equation}
The convergence of $N(t,\omega)$ to $-e_3$
as $t \rightarrow + \infty$ is equivalent to the convergence of
$M_2(t,\omega)$ to $M_f(\omega)$ as $t \rightarrow + \infty$.

\subsection{Lyapunov feedback}

Let us consider the following Lyapunov-like functional
\begin{multline}~\label{lyap:eq}
 \mathcal{L}(N) :=   \frac{\|N+e_3\|_{H^1}^2}{2} \\=
\int\limits_{\wmin}^{\wmax} \Big( \frac{1}{2} \Big\| \frac{\partial
N}{\partial \omega} \Big\|^2 + 1 + \langle N , e_3 \rangle \Big)
d\omega.
\end{multline}
The function $\mathcal{L}$ is defined for any $N \in
H^1((\omega_*,\omega^*),\SSS^2)$ and takes its minimal value on this
space at the point $N=-e_3$ with $\mathcal{L}(-e_3)=0$.
For any solution of (\ref{equ:N}), some  computations show that
$$\frac{d\mathcal{L}}{dt}[N(t)] = \sum_{i=1}^2 u_i(t) H_i[t,N(t)]$$
where, for $i=1,2$ one has
\begin{multline*}
H_i[t,N]:=
\int_{\omega_*}^{\omega^*} \bigg[ \left\langle
\frac{dN}{d\omega}(\omega) , \Big( \frac{\partial F}{\partial
\omega}(t,\omega) e_i \Big) \wedge N(\omega) \right\rangle
\\+
  \left\langle e_3 , \Big( F(t,\omega) e_1 \Big) \wedge
N(\omega) \right\rangle \bigg] d\omega.
 \end{multline*}
 Hence, with the feedback laws
\begin{equation} \label{def:u(N)}
u_i(t,N):=-H_i[t,N], \forall i\in\{1,2\},
\end{equation}
it follows that
\begin{equation} \label{uvL2}
\frac{d\mathcal{L}}{dt}[N(t)]=-u_1(t,N)^2 - u_2(t,N)^2 \leqslant 0.
\end{equation}
As in \cite{beauchard-et-al:auto11}, we have the following result.
\begin{Prop} \label{Solutions}
For every initial condition $N_0 \in H^1((\wmin, \wmax), \SSS^2)$,
the closed loop system (\ref{equ:N}), (\ref{def:u(N)}) has a unique
solution $N \in C^1 \left([0, \infty), H^1\left( (\wmin, \wmax),
\mathbb{R}^3 \right) \right)$ such that $N(0)=N_0$.
\end{Prop}

\subsection{Closed-loop simulations}

We assume here $\wmin=0$, $\wmax=1$ and we solve numerically the
$T$-periodic system (1) with  the feedback law $(\tilde{u}_1,
\tilde{u}_2)$ given by  (\ref{uv:eq}), (\ref{def:u(N)}).
The closed-loop  simulation is performed for $t\in[0, T_f]$, $T_f=20 T$ and $T=2\pi/(\wmax-\wmin)$. The
$\omega$-profile $[\wmin,\wmax]\ni \omega \mapsto
(x(t,\omega),y(t,\omega),z(t,\omega))$ is discretized $\{1,\ldots,
N+1\}\ni k \mapsto (x_k(t),y_k(t),z_k(t))$ with a regular mesh of
step $\epsilon_N = \frac{\wmax-\wmin}{N}$ with $N=100$. In other
words,  one has a set of discrete values $\{\omega_i, i=1, \ldots,
N+1\}$, where $\omega_i =\wmin + (i-1) \epsilon_N$.

 We have
checked  that the closed-loop simulations are almost identical for
$N=100$ and $N=200$.  In the feedback law (16), the integral versus
$\omega$ is computed assuming that $(x,y,z)$ and
$(x^\prime,y^\prime,z^\prime)$ are constant over
   { $](k-\tfrac{1}{2})\epsilon_N,(k+\tfrac{1}{2})\epsilon_N[$},
   their values being $(x_k,y_k,z_k)$ and { $\left(\tfrac{x_{k+1}-x_{k-1}}{2\epsilon_N},
   \tfrac{y_{k+1}-y_{k-1}}{2\epsilon_N},\frac{z_{k+1}-z_{k-1}}{2\epsilon_N}\right)$.}
The obtained differential  system is  of  dimension $3(N+1)$. It is
integrated via an explicit Euler scheme with a step size $h=
T/1000$. We have tested that $h=T/2000$ yields  almost the same
numerical solution at $t=T_f=20T$. After each time-step the new
values of $(x_k,y_k,z_k)$ are  normalized to remain in $\SSS^2$. The
initial $\omega$-profile $M_0(\omega)$ of $(x,y,z)\in\SSS^2$ is
given by $x_0 =0, y_0=-\sqrt{1-z_0^2}$, where
$z_0=-\cos(\frac{\pi}{8})+0.05\left(1-\cos(\frac{\pi}{8})\cos(\omega
\frac{\pi}{2})\right)$. The desired final profile $M_f(\omega)$ is
given by $x_f=-\sqrt{1-z_f^2}, y_f=0$, where
$z_f=-\cos(\frac{\pi}{16})+0.1\left(1-\cos(\frac{\pi}{16})\sin(\omega
\frac{\pi}{4}) \right)$.

The map $R(\omega)$ is constructed for the discrete set $\{\omega_i,
i=1, \ldots, N+1\}$, in the following way. For $i=1$, one takes
$r_3(\omega_1) = M_0(\omega_1)$. Now choose a vector $\theta$ among
the vectors of the canonical basis in a way that $\langle \theta,
M_f(\omega_1) \rangle$ is the minimum value. Construct
$r_2(\omega_1) = \frac{1}{\| \theta \wedge r_3(\omega_1) \|}
\left(\theta \wedge r_3(\omega_1) \right)$. Then one may take
$r_1(\omega_2) = r_2(\omega_2) \wedge r_3(\omega_2)$. Now, for $i=2,
3, \ldots, N+1$ one chooses $r_3(\omega_i) = M_0(\omega_i)$, $\theta
= -r_1(\omega_{i-1})$, $r_2(\omega_i) = \frac{1}{\| \theta \wedge
r_3(\omega_i)  \|} \left( \theta \wedge r_3(\omega_i) \right)$ and
$r_1(\omega_i) = r_2(\omega_i)  \wedge r_3(\omega_i) $, and so on.
The orthogonal matrix $\bar R(\omega)$ formed by the column
vectors $r_1$, $r_2$, $r_3$ is then transposed to obtain
$R(\omega)$.

Figures~\ref{fig:LyapControl} and~\ref{fig:BeginEnd} summarize
the main convergence issues for these
choices of initial profile $M_0$ and of the desired final profile
$M_f$. The convergence speed is rapid at the beginning and tends to
decrease at the end. We start with ${\mathcal L}(0)\approx 0.1929$.
We get ${\mathcal L}(20T)\approx 0.0032$. This numerically observed
convergence is confirmed by Theorem~\ref{MainThm} here below.

\begin{figure}[bht]
  \centerline{\includegraphics[width=\columnwidth]{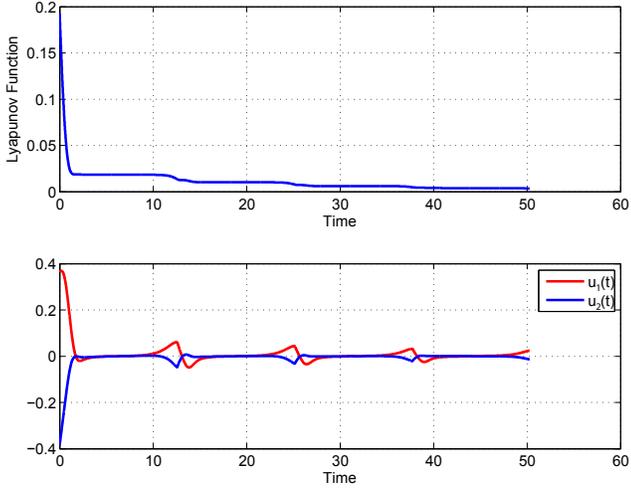}}
  \caption{Lyapunov function $\mathcal{L}(t)$ defined by (\ref{lyap:eq})
  and  the closed-loop control $({u}_1,{u}_2)$ defined by~(\ref{def:u(N)})}\label{fig:LyapControl}
\end{figure}

\begin{figure}[bht]
  \centerline{\includegraphics[width=\columnwidth]{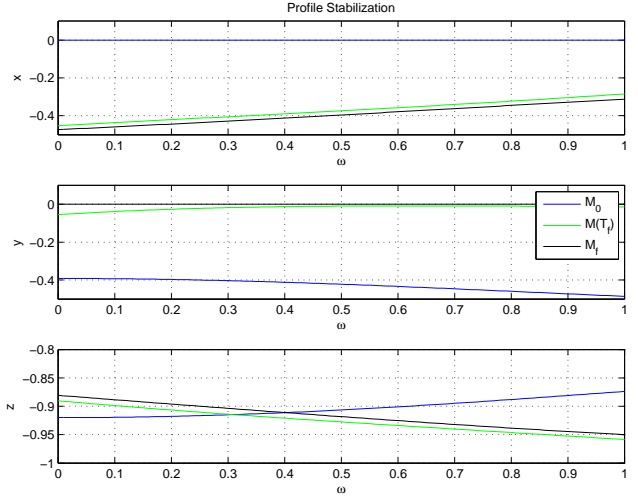}}
  \caption{Initial ($t=0$) and final ($t=T_f$) $\omega$-profiles for
  $x$, $y$ and $z$ solutions of the closed-loop system \ref{dyn:eq1}
  with the  feedback (\ref{uv:eq}), (\ref{def:u(N)}).}\label{fig:BeginEnd}
\end{figure}


\section{Main Result} \label{main:sec}

\subsection{Local stabilization}

\begin{thm} \label{MainThm}
For every $M_f \in H^1((\omega_*,\omega^*),\SSS^2)$ with
\begin{equation} \label{equateur}
\langle M_f(\omega), e_3 \rangle \neq 0, \forall \omega \in [\omega_*,\omega^*],
\end{equation}
there exists $\delta_1>0$ such that, for every $N_0 \in H^1((\omega_*,\omega^*),\SSS^2)$
with $\|N_0 + e_3\|_{H^1} \leq \delta_1$, the solution of the closed loop system (\ref{equ:N}),
(\ref{def:u(N)}) with initial condition $N(0,\omega)=N_0(\omega)$
satisfies $N(t) \rightharpoonup -e_3$ weakly in $H^{1}(\omega_*,\omega^*)$
when  $t \rightarrow + \infty$.
\end{thm}

The above theorem has the following corollary.

\begin{Cor}
For every $M_f \in H^1((\omega_*,\omega^*),\SSS^2)$ with (\ref{equateur}),
there exists $\delta_2>0$ such that,
for every $M_0 \in H^1((\omega_*,\omega^*),\SSS^2)$ with
$\|M_0-M_f\|_{H^1}<\delta_2$, the solution of the system
(\ref{dyn:eq1}) with the initial condition $M(0,\omega)=M_0(\omega)$
and the feedback law given by (\ref{uv:eq}), (\ref{def:u(N)})
satisfies $M((2kT)^+) \rightharpoonup M_f$ weakly in $H^{1}$ when $k \rightarrow + \infty$.
In particular,
$$\lim\limits_{k \rightarrow + \infty} \| M((2kT)^+,.) - M_f \|_{L^\infty(\omega_*,\omega^*)} =0.$$
\end{Cor}

The remaining part of this section is devoted to the proof of Theorem~\ref{MainThm}

\subsection{LaSalle invariant set}

The first step of our proof consists in checking that, locally, the
invariant set is reduced to $\{-e_3\}$.

\begin{Prop} \label{Prop:Inv}
For every $M_f \in H^1((\omega_*,\omega^*),\SSS^2)$ with (\ref{equateur}),
there exists $\delta >0$ such that, for every $N_0 \in H^1((\omega_*,\omega^*),\SSS^2)$ with
$\|N_0+e_3\|<\delta$, the map $t \mapsto \mathcal{L}[N(t)]$ is constant on $[0,+\infty)$ if and only if $N_0=-e_3$.
\end{Prop}

\noindent \textbf{Proof:} Let us assume that $t \mapsto \mathcal{L}[N(t)]$ is constant.
Then $u_1=u_2=0$ and $N(t,\omega) \equiv N_0(\omega)$ (see (\ref{uvL2}) and (\ref{equ:N})).
Thus, for every $j \in \{1,2\}$ and $t \in [0,+\infty)$
\begin{equation} \label{inv1}
\begin{array}{rl}
0 = & \int_{\omega_*}^{\omega^*} \left[ \left\langle N_0'(\omega)
, \Big( \frac{\partial F}{\partial \omega}(t,\omega) e_j \Big)
\wedge N_0(\omega) \right\rangle  \right.
\\  &  \left. +
\left\langle e_3 , \Big( F(t,\omega) e_j \Big) \wedge N_0(\omega)
\right\rangle \right] d\omega.
 \end{array}
\end{equation}
For $t \in [0,T]$, $\sigma(t)=t$ so
$F(t,\omega)=\sum_{k=0}^{\infty} \frac{t^k \omega^k}{k!} R(\omega)S^k$
and
$\frac{\partial F}{\partial \omega}(t,\omega)=
\sum\limits_{k=0}^{\infty} \frac{t^k \omega^k}{k!} R'(\omega)S^k +
\sum\limits_{k=1}^{\infty} \frac{t^{k} \omega^{k-1}}{(k-1)!}
R(\omega) S^k$.
Developing (\ref{inv1}) in power series expansions of $t$
and using (\ref{def:S}), we obtain, for every $j \in \{1,2\}$ and $k \geqslant 1$,
  \[
 \begin{array}{c}
\int\limits_{\omega_*}^{\omega^*} \left\langle N_0'(\omega) , \left[
\left( \frac{\omega^k}{k!} R'(\omega) + \frac{\omega^{k-1}}{(k-1)!}
R(\omega) \right)  e_j \right]
 \wedge N_0(\omega) \right\rangle \\
+ \left\langle e_3 , \Big( \frac{\omega^k}{k!} R(\omega)  e_j \Big)
\wedge N_0(\omega) \right\rangle
 d\omega=0.
 \end{array}
  \]
By linearity, the following equality holds, for every $Q \in
\mathbb{R}[X]$ and $j \in \{1,2\}$
\begin{equation} \label{inv2}
 \begin{array}{c}
\int\limits_{\omega_*}^{\omega^*}
\left\langle N_0'(\omega) , \Big[
\left( Q(\omega) R'(\omega) + Q'(\omega) R(\omega) \right)  e_j
\Big]
 \wedge N_0(\omega) \right\rangle \\
+ \left\langle e_3 , \Big[ Q(\omega) R(\omega)  e_j \Big] \wedge
N_0(\omega) \right\rangle
 d\omega=0.
  \end{array}
\end{equation}
Thanks to the density of polynomial functions in
$H^1((\omega_*,\omega^*),\mathbb{C})$, the previous equality holds
for every $Q \in H^1((\omega_*,\omega^*),\mathbb{R})$.
Let us recall the relations
$\langle X , Y \wedge Z \rangle = \langle Y , Z \wedge X \rangle$
and
$\langle MX,Y \rangle = \langle X , M^\top Y\rangle$, $\forall X,Y,Z \in \mathbb{R}^3$
and $M \in \mathcal{M}_3(\mathbb{R})$, where $M^\top$ denotes the transposed matrix of $M$.
Then, the equality (\ref{inv2}) may also we written
\begin{equation} \label{inv2_bis}
 \begin{array}{c}
\int\limits_{\omega_*}^{\omega^*}
\left\langle Q(\omega) e_j , R'(\omega)^\top [ N_0(\omega) \wedge N_0'(\omega)] \right\rangle \\
+\left\langle Q'(\omega) e_j , R(\omega)^\top [ N_0(\omega) \wedge N_0'(\omega)] \right\rangle \\
+ \left\langle Q(\omega) e_j , R(\omega)^\top [N_0(\omega) \wedge e_3 ]\right\rangle
 d\omega=0.
  \end{array}
\end{equation}
for every $Q \in H^1((\omega_*,\omega^*),\mathbb{R})$ and $j \in \{1,2\}$.
By linearity, we deduce that
\begin{equation} \label{inv3}
 \begin{array}{c}
\int\limits_{\omega_*}^{\omega^*}
\left\langle \mathcal{Q}(\omega) ,
R'(\omega)^\top [N_0(\omega) \wedge N_0'(\omega)] + R(\omega)^\top [N_0(\omega) \wedge e_3]
\right\rangle
\\
+ \left\langle \mathcal{Q}'(\omega) , R(\omega)^\top [ N_0(\omega) \wedge N_0'(\omega) ]
\right\rangle d\omega = 0
  \end{array}
\end{equation}
for every $\mathcal{Q} \in H^1((\omega_*,\omega^*),\mathbb{V})$ where $\mathbb{V}:=\text{Span}(e_1,e_2)$.
Let $\mathbb{P}:\mathbb{R}^3 \rightarrow \mathbb{V}$ be the orthogonal projection on $\mathbb{V}$.
The previous equality is equivalent to
\begin{equation} \label{syst_inv_1}
\left\lbrace \begin{array}{l}
\mathbb{P} R(\omega)^\top [ N_0(\omega)\wedge e_3 - N_0(\omega)\wedge N_0''(\omega) ] =0
\text{ in } H^{-1} \\
R(\omega)^\top[ N_0(\omega)\wedge N_0'(\omega)]=0  \text{ at } \omega=\omega_* \text{ and } \omega^*.
\end{array}\right.
\end{equation}
Here, $H^{-1}$ denotes the dual space of $H^1_0(\omega_*,\omega^*)$ for the $L^2$-scalar product;
the first equation has to be understood in the distribution sens.
Thanks to $\|N_0(\omega)\| \equiv 1$, we have $\| R(\omega)^\top[ N_0(\omega)\wedge N_0'(\omega)] \| \equiv \|N_0'(\omega)\|$.
Thus, the second line of (\ref{syst_inv_1}) is equivalent to $N_0'=0$ at $\omega_*$ and $\omega^*$.
Notice that $\mathbb{P}R(\omega)^\top |_\mathbb{V}$ is bijective on $\mathbb{V}$ for every $\omega \in [\omega_*,\omega^*]$.
Indeed, thanks to (\ref{def:R}) and (\ref{equateur}), we have
$$\begin{array}{ll}
\text{Range}[\mathbb{P}R(\omega)^\top|_\mathbb{V} ]^{\perp}
& = \text{Ker}[\mathbb{P}R(\omega)|_\mathbb{V}] \\
& = \{ v \in \mathbb{V} ; R(\omega) \in \mathbb{R} e_3 \}\\
& = \mathbb{V} \cap \mathbb{R} M_f(\omega) = \{0\}.
\end{array}$$
Moreover, $(\mathbb{P}R(\omega)^\top|_\mathbb{V})^{-1} \in H^1$, thus (\ref{syst_inv_1}) gives
\begin{equation} \label{syst_inv_2}
\left\lbrace \begin{array}{l}
- N_0''\wedge e_3 + N_0 \wedge e_3  = g \text{ in } H^{-1}((\omega_*,\omega^*),\mathbb{V}),\\
N_0'\wedge e_3=0  \text{ at } \omega_*,\omega^*,
\end{array}\right.
\end{equation}
where
$$g(\omega):=- (\mathbb{P}R(\omega)^\top|_\mathbb{V})^{-1} \mathbb{P}R(\omega)^\top[ N_0''(\omega) \wedge (N_0(\omega)+e_3) ].$$
Therefore, there exists $C_1=C_1(\omega_*,\omega^*)>0$ such that
\begin{equation} \label{reg_ellipt}
\| N_0(\omega)\wedge e_3 \|_{H^1} \leqslant C_1 \| g \|_{H^{-1}}.
\end{equation}
Thanks to (\ref{R_borne}), there exists $C_2=C_2(\omega_*,\omega^*,\|M_f\|_{H^1})>0$ such that
$$\| g \|_{H^{-1}} \leqslant C_2 \| N_0 + e_3 \|_{H^1}^2.$$
When $N_0$ is close enough to $-e_3$ in $H^1$, then $\| N_0 \wedge e_3 \|_{H^1}$ and
$\| N_0 + e_3 \|_{H^1}$ are equivalent norms and then (\ref{reg_ellipt}) gives
$$\| N_0(\omega)\wedge e_3 \|_{H^1} \leqslant C_3 \| N_0(\omega) \wedge e_3 \|_{H^1}^2$$
for some constant $C_3=C_3(\omega_*,\omega^*,\|M_f\|_{H^1})>0$.
This implies $N_0 \wedge e_3 = 0$, i.e. $N_0=-e_3$. $\hfill \Box$

\begin{rk}
For $M_f \equiv e_1$, any constant function $N_0$ with values in $\text{Span}(e_2,e_3)$
belongs to the invariant set (see (\ref{syst_inv_1})). Thus, an assumption of the type (\ref{equateur})
is required for our strategy to work.
\end{rk}

\subsection{Convergence proof}

For the proof of Theorem \ref{MainThm}, we need the following result.

\begin{Prop} \label{Prop:continuiteH1faible}
Take  $M_f \in
H^1((\omega_*,\omega^*),\SSS^2)$ and
$R\in
H^1((\omega_*,\omega^*),SO(3))$ as in Proposition \ref{Def:R}. Let $(N_n^0)_{n \in \mathbb{N}}$ a sequence of
$H^1((\omega_*,\omega^*),\SSS^2)$ and $N_\infty^0 \in
H^1((\omega_*,\omega^*),\SSS^2)$ such that $N_n^0 \rightharpoonup
N_\infty^0$ weakly in $H^1$ and $N_n^0 \rightarrow N_\infty^0$
strongly in $L^2$. Let $\alpha \in [0,2T]$ and $(\tau_n)_{n \in
\mathbb{N}}$ be a sequence of $[0,2T)$ such that $\tau_n \rightarrow
\alpha$. Let $N_n$ (resp. $N_\infty$) be the solutions of the closed
loop system (\ref{equ:N}), (\ref{def:u(N)}) associated to the
initial condition $N_n(\tau_n)=N_n^0$ (resp.
$N_\infty(\alpha)=N_\infty^0$). Then, we have
$N_n(t) \rightharpoonup N_\infty(t)$ weakly in $H^1$, $\forall t > \alpha$, and
$u_j[t,N_n(t)] \rightarrow u_j[t,N_\infty(t)]$, $\forall t > \alpha, \forall j \in \{1,2\}$.
\end{Prop}

\noindent \textbf{Proof:} The sequence $(N_n^0)_{n \in \mathbb{N}}$ is bounded in $H^1$ and
$\mathcal{L}[N_n(t)] \leqslant \mathcal{L}[N_n^0]$, for every $t \in [\tau_n,+\infty)$ and $n \in \mathbb{N}$
so there exists $\mathcal{M}_0>0$ such that
$\|N_n(t)\|_{H^1} \leqslant \mathcal{M}_0$, for every $t \in [\tau_n,+\infty)$ and $n \in \mathbb{N}$.
The function $t \in \mathbb{R} \mapsto F(t,.) \in H^1((\omega_*,\omega^*),\mathcal{M}_3(\mathbb{R}))$
defined by (\ref{def:F}) is continuous and $2T$-periodic, thus, there exists $\mathcal{M}_1>0$ such that
$\|F(t,.)\|_{H^1} \leqslant \mathcal{M}_1$, for every $t \in \mathbb{R}$.
Thanks to (\ref{def:u(N)}), we have
$|u_j[t,N_n(t)]| \leqslant 2 \mathcal{M}_1 \mathcal{M}_0$, for every $n \in \mathbb{N}$ and $t \in [\tau_n,+\infty)$.
We deduce from (\ref{equ:N}) that
$\Big\| \frac{\partial N_n}{\partial t}(t) \Big\|_{L^2} \leqslant 4 \mathcal{M}_1 \mathcal{M}_0$
for every $t \in [\tau_n,+\infty)$ and $n \in \mathbb{N}$.
The end of the proof is as in \cite{beauchard-et-al:auto11}. \hfill $\Box$

\noindent \textbf{Proof of Theorem \ref{MainThm}:}
The proof is as in \cite{beauchard-et-al:auto11}.
One may replace Barbalat's lemma by the Lebesgue reciprocal theorem, in the following way.
Thanks to (\ref{uvL2}), $t \mapsto u_i[t,N(t)]$ belongs to $L^2(0,+\infty)$,
thus, for any diverging sequence of times $k_n$,
the sequence $(t \in (0,+\infty) \mapsto u_i[2k_nT+t,N(2k_nT+t)])_{n \in \mathbb{N}}$
converges to zero in $L^2(0,+\infty)$.
Therefore, there exists a subset $\mathcal{N} \subset (0,+\infty)$ with zero Lebesgue measure such that
$u_j[t,N(2k_nT+t) \rightarrow 0$ for every
$t \in (0,+\infty)-\mathcal{N}$ and $j \in \{1,2\}$. \hfill $\Box$

\section{Concluding remark}
Open-loop "impulse-train" control are combined  with Lyapunov feedback to steer an initial profile $[\omega_*,\omega^*]\ni\omega\mapsto M(0,\omega)$ of the Bloch-sphere system~\eqref{dyn:eq1} towards an arbitrary target profile $[\omega_*,\omega^*]\ni\omega\mapsto M_f(\omega)$. Convergence is proved to be local for any target profile belonging  either to  the south or to  the north hemisphere. We guess that our convergence proof could be extended to the case where $M_f$ intersects transversely  the equator and thus where $M_f$  is not  confined in only one hemisphere.

\begin{ack}                               
KB and PR were partially supported by the ``Agence Nationale de la
Recherche'' (ANR), Projet Blanc C-QUID number BLAN-3-139579 and
Project Blanc EMAQS number ANR-2011-BS01-017-01. PSPS is
Partially supported by CNPq -- Conselho Nacional de Desenvolvimento
Cientifico e Tecnologico, and FAPESP-Funda\c{c}\~ao de Amparo a
Pesquisa do Estado de S\~ao Paulo.
\end{ack}


\end{document}